
\input amstex
\magnification = 1200
\documentstyle{amsppt}
\NoBlackBoxes

 \def\today{\ifcase\month\or
  January\or February\or March\or April\or May\or June\or
  July\or August\or September\or October\or November\or December\fi
  \space\number\day, \number\year}

\document

\topmatter
\title On representations of integers by  indefinite ternary
       quadratic forms\endtitle
\rightheadtext{Representations by  ternary quadratic forms}
\author Mikhail  Borovoi \endauthor
\address  Raymond and Beverly Sackler School of Mathematical Sciences,
Tel Aviv University, 69978 Tel Aviv, Israel \endaddress
\email borovoi\@math.tau.ac.il\endemail
\thanks Partially supported by the Hermann Minkowski Center for Geometry
\endthanks
\keywords Ternary quadratic forms \endkeywords
\subjclass Primary 11E12 \endsubjclass
\abstract
Let $f$ be an indefinite ternary integral quadratic form and let $q$ be
a nonzero integer such that $-q\text{det}(f)$ is not a square.
  Let $N(T,f,q)$ denote the number of integral solutions 
of the equation $f(x)=q$ where $x$ lies in the ball of radius $T$
centered at the origin. We are interested in the 
asymptotic behavior of $N(T,f,q)$ as $T\to\infty$.
  We deduce from the results of our joint paper with Z.Rudnick that 
$N(T,f,q)\sim cE_{HL}(T,f,q)$ as $T\to\infty$, where $E_{HL}(T,f,q)$
is the Hardy-Littlewood expectation (the product of local densities) and 
$0 \le c \le 2$. We give examples of $f$ and $q$ such that $c$ takes  
values 0,1,2.  
\endabstract

\endtopmatter

  \def\bet{{\beta}} 
 \def\del{{\delta}}
   
\def\kap{{\kappa}}                   
\def\lam{{\lambda}}    \def\Lam{{\Lambda}}   
    \def\om{{\omega}}

\def\R{\bold {R}}  \def\Q{\bold {Q}}  \def\F{{\bold F}}
 \def\Fp{\bold F_p}
\def\Z{\bold {Z}} \def\A{{\bold A}}
   
\def\Zhat{\hat{\Z}}

\def\cU{{\Cal U}}
\def\o{{\Cal O}}
\def\oA{{\o_\A}}
\def\ov{{\o_v}}

\def\SS{\frak S}

\def\SO{\text{SO}}

\def\textrm#1{\text{\rm #1}}

            \def\kb{{\bar k}}

 \def\o{{\Cal O}}

\def\Gal{\textrm{Gal}}       
           
  \def\GL{\textrm{GL}}

\def\tors{_{\textrm{tors}}}

\def\Spin{\textrm{Spin}}
\def\Hom{\textrm{Hom}}
\def\Vol{\textrm{Vol}}

\def\half{{1\over 2}}

\heading 0. Introduction \endheading

Let $f$ be a nondegenerate indefinite integral-matrix quadratic form of $n$
variables:
$$
f(x_1,\dots,x_n)=\sum_{i=1, j=1}^n a_{ij}x_ix_j,\quad a_{ij}\in\Z,\quad 
a_{ij}=a_{ji}\,.
$$
  Let $q\in \Z$, $q\neq 0$.
  Let $W=\Q^n$.
  Consider the affine quadric $X$ in $W$ defined by the equation
$$
f(x_1,\dots,x_n)=q\,.
$$
  We wish to count the representations of $q$ by the quadratic form $f$,
that is the integer points of $X$.

Since $f$ is indefinite, the set $X(\Z)$ can be infinite. 
  We fix a Euclidean norm $|\cdot|$ on $\R^n$. 
  Consider the counting function
$$
N(T,X)=\#\{x\in X(\Z): |x|\leq T\}
$$
where $T\in \R,\ T>0$.
  We are interested in the asymptotic behavior of $N(T,X)$ as $T\to\infty$.

When $n\ge 4$, the counting function $N(T,X)$ can be approximated by the 
product of local densities.
  For a prime $p$ set
$$
\mu_p(X)=\lim_{k\to\infty}\frac{\#X(\Z/p^k\Z)}{(p^k)^{n-1}}\;.
$$
  For almost all $p$ it suffices to take $k=1$:
$$
\mu_p=\frac{\#X(\F_p)}{p^{n-1}}\;.
$$
Set $\SS(X)=\prod_p\mu_p(X)$, this product converges absolutely 
(for $n\ge 4$), it is called the singular series.
Set
$$
\mu_\infty(T ,X)=\lim_{\varepsilon\to 0}
\frac{\Vol \{x\in \R^n : |x|\leq T,\ |f(x)-q|<\varepsilon/2,\}}
{\varepsilon} \;,
$$
it is called the singular integral.

\proclaim{Theorem} For $n\ge 4$
$$
N(T,X)\sim \SS(X)\mu_\infty(T,X) \text{ as } T\to\infty.
$$
\endproclaim

This theorem follows from results of [BR], 6.4 (based on analytical
results of [DRS], [EM], [EMS]). For certain non-Euclidean norms it
was earlier proved by the Hardy-Littlewood circle method, cf. [Da], [Est].

We are interested here in the case $n=3$, a ternary quadratic form.
  This case is beyond the range of the Hardy-Littlewood circle method.
  Set $D=\det(a_{ij})$. 
We assume that $-qD$ is not a square.
  Then the product $\SS(X)=\prod\mu_p(X)$ conditionally converges 
(see Sect.~1
below), but in general $N(T,X)$ is not asymptotically 
$\SS(X)\mu_\infty(T,X)$.
   From results of [BR] it follows that 
$$
N(T,X)\sim c_X \SS(X)\mu_\infty(T,X) \text{ as } T\to\infty
$$
with $0\leq c_X\leq 2$, see details in  Subsection 1.5 below. 
  We wish to know what values can take $c_X$.

A case when $c_X=0$ was already known to Siegel, see also [BR], 6.4.1.
  Consider the quadratic form 
$$
f_1(x_1,x_2,x_3)=-9x_1^2+2x_1x_2+7x_2^2+2x_3^2\;,
$$
and take $q=1$. Let $X$ be defined by $f_1(x)=q$. 
  Then $f_1$ does not represent $1$ over $\Z$, so $N(T,X)=0$ for all $T$.
  On the other hand, $f_1$ represents $1$ over $\R$ and over $\Z_p$ for
all $p$, and $\SS(X)\mu_\infty(T,X) \to \infty$ as $T\to \infty$.
Thus $c_{X}=0$ (see details in Sect.~2).

We show that $c_X$ can take value 2.
  Recall that two integral quadratic forms $f,f'$ are in the same genus, if
they are equivalent over $\R$ and over $\Z_p$ for every prime $p$, cf. e.g.
[Ca].

\proclaim{Theorem 0.1} Let $f$ be an indefinite integral-matrix
ternary quadratic
form, $q\in\Z$, $q\neq 0$, and let $X$ be the affine quadric defined by the 
equation $f(x)=q$.
  Assume that $f$ represents $q$ over $\Z$ and there exists a quadratic form
$f'$ in the genus of $f$, such that $f'$ does not represent $q$ over $\Z$.
Then $c_X=2$:
$$
N(T,X)\sim 2 \SS(X)\mu_\infty(T,X) \text{ as } T\to\infty.
$$
\endproclaim
Theorem 0.1 will be proved in Sect. 3.

\proclaim{\it Example 0.1.1} \rm Let  
$f_2(x_1,x_2,x_3)=-x_1^2+64x_2^2+2x_3^2,\ q=1$.
  Then $f_2$ represents 1 ($f_2(1,0,1)=1$) and the quadratic form $f_1$ 
considered above is in the genus of $f_2$ (cf. [CS], 15.6).
 The form $f_1$ does not represent 1.
Take $|x|=(x_1^2+64x_2^2+2x_3^2)^{1/2}$.
  By Theorem 0.1 $c_{X}=2$ for the variety $X: f_2(x)=1$.
  Analytic and numeric calculations give 
$2\SS(X)\mu_\infty(T,X)\sim 0.794 T$.
  On the other hand, numeric calculations give for $T=10,000$ 
$N(T,X)/T=0.8024$.
\endproclaim

\comment
\proclaim{Remark 0.1.2} \rm In Example 0.1.1 the equality $c_X=2$
follows from Siegel's weight formula 
(and the main result of [DRS], [EM], [EMS]), 
because there are only two classes in the genus of $f_2$
and the other class does not represent $q$. 
  However when in Theorem 0.1 the genus of $f$ contains more than 2
classes, the assertion of Theorem 0.1 seems not to follow immediately
from Siegel's weight formula.
\endproclaim
\endcomment

We also show that $c_X$ can take the value 1.

\proclaim{Theorem 0.2} Let $f$ be an indefinite integral-matrix 
ternary quadratic
form, $q\in\Z$, $q\neq 0$, and let $X$ be the affine quadric defined by the 
equation $f(x)=q$.
  Assume that $X(\R)$ is two-sheeted (has two connected components). 
Then $c_X=1$:
$$
N(T,X)\sim \SS(X)\mu_\infty(T,X) \text{ as } T\to\infty.
$$
\endproclaim
Theorem 0.2 will be proved in Sect. 4.

\proclaim{\it Example 0.2.1} \rm Let $f_2$ and $|x|$ be as in Example 0.1.1, 
$q=-1$, $X:\ f_2(x)=q$.
Then $X(\R)$ has two connected components, and by Theorem 0.2 $c_X=1$.
  Analytic and  numeric calculations give 
$\SS(X)\mu_\infty(T,X)\sim 0.7065 T$. 
  On the other hand, numeric calculations give for 
$T=10,000\ N(T,X)/T=0.7048$.
\endproclaim

\proclaim{Question 0.3}\rm Can $c_X$ take values other than 0,\ 1,\ 2? 
\endproclaim

\proclaim{Remark 0.4} 
\rm It seems that Theorems 0.1 and 0.2 also can be proved using a result
of Kneser ([Kn], Satz 2) together with Siegel's weight formula [Si] and the 
results of [DRS], [EM], [EMS].
\endproclaim

The plan of the paper is the following. 
  In Section 1 we expose results
of [BR] in the case of 2-dimensional affine quadrics. 
  In Section 2 we treat in detail the example of $c_X=0$.
  In Section 3 we prove Theorem 0.1.
  In Section 4 we prove Theorem 0.2.

{\it Acknowledgements.} This paper was partly written when the author
was visiting Sonderforschungsbereich 343 ``Diskrete Strukturen in der
Mathematik'' at Bielefeld University, and  I am grateful to SFB 343 
for hospitality and support. 
  I thank Rainer Schulze-Pillot and John S. Hsia for useful e-mail 
correspondence.
  I am  grateful to Ze\'ev Rudnick  for useful discussions and
help in analytic calculations.

\head 1. Results of [BR] in the case of ternary quadratic forms
                    \endhead

\subhead 1.0 \endsubhead
  Let $f$ be an indefinite ternary integral-matrix quadratic form
$$
f(x_1,x_2,x_3)=\sum_{i,j=1}^3 a_{ij}x_ix_j,\quad a_{ij}\in\Z,\quad 
a_{ij}=a_{ji}\,.
$$
  Let $q\in \Z$, $q\neq 0$. Let $D=\det(a_{ij})$. 
We assume that $-qD$ is not a square.

Let $W=\Q^3$ and let $X$ denote the affine variety in $W$ defined by
the equation $f(x)=q$, where $x=(x_1,x_2,x_3)$.
  We assume that $X$ has a $\Q$-point $x^0$. 
  Set $G=\Spin(W,f)$, the spinor group of $f$.
  Then $G$ acts on $W$ on the left, and $X$ is an orbit (a homogeneous
space) of $G$.

\subhead 1.1. Rational points in adelic orbits \endsubhead
Let $\A$ denote the ad\`ele ring of $\Q$. 
  The group $G(\A)$ acts on $X(\A)$; let $\oA$ be an orbit.
  We are interested whether $\oA$ has a $\Q$-rational point.

Let $W'$ denote the orthogonal complement of $x^0$ in $W$, 
and let  $f'$ denote the restriction of $f$ to $W'$.
  Let $H$ be the stabilizer of $x^0$ in $G$, then
$H=\Spin(W',f')$.
  Since $\dim W'=2$, the group $H$ is a one-dimensional torus.

We have $\det f' ={D}/{q}$, so up to multiplication by a square
$\det f'=qD$.
  It follows that up to multiplication by a scalar, $f'$ is equivalent
to the quadratic form $u^2+qDv^2$. 
  Set $K=\Q(\sqrt{-qD})$, then $K$ is a quadratic extension of $\Q$,
because $-qD$ is not a square.
  The torus $H$ is anisotropic over $\Q$ (because $-qD$ is not a square),
and $H$ splits over $K$.
  Let $X_*(H_K)$ denote the cocharacter group of $H_K$,
$X_*(H_K)=\Hom(\Bbb G_{m,K},H_K)$; then $X_*(H_K)\simeq\Z$.
  The non-neutral element of $\Gal(K/\Q)$ acts
on $X_*(H_K)$ by multiplication by $-1$.

Let $\oA$ be an orbit of $G(\A)$ in $X(\A)$, $\oA=\prod \ov$ where $\ov$
is an orbit of $G(\Q_v)$ in $X(\Q_v)$, $v$ runs over the places of $\Q$, 
and $\Q_v$ denotes the completion of $\Q$ at $v$.
  We define local invariants $\nu_v(\ov)=\pm 1$.
  If $\ov=G(\Q_v)\cdot x^0$, then we set $\nu_v(\ov)=+1$, if not, we set
$\nu_v(\ov)=-1$. 
  Then $\nu_v(\ov)=+1$ for almost all $v$.
  We define $\nu(\oA)=\prod\nu_v(\ov)$ where $\oA=\prod \ov$.
  Note that the local invariants $\nu_v(\ov)$ depend on the choice of the
rational point $x^0\in X(\Q)$; one can prove, however, that their product
$\nu(\oA)$ does not depend on $x^0$.

Let $x\in X(\A)$.
  We set $\nu(x)=\nu(G(\A)\cdot x)$.
  Then $\nu(x)$ takes values $\pm 1$; it is a locally constant function on 
$X(\A)$, because the orbits of $G(\A)$ are open in $X(\A)$.

For $x\in X(\A)$ define $\del(x)=\nu(x)+1$.
In other words, if $\nu(x)=-1$ then $\del(x)=0$, and if $\nu(x)=+1$
then $\del(x)=2$.
  Then $\del$ is a locally constant function on $X(\A)$.

\proclaim{Theorem 1.1.1}
  An orbit $\oA$ of $G(\A)$ in $X(\A)$ has a  $\Q$-rational point if 
and only if $\nu(\oA)=+1$.
\endproclaim

Below we will deduce Theorem 1.1.1 from [BR], Thm. 3.6.

\subhead 1.2. Proof of Theorem 1.1.1 \endsubhead
 For a torus $T$ over a field $k$ of characteristic 0 we define a finite
abelian group $C(T)$ as follows:
$$
C(T)=(X_*(T_\kb)_{\Gal(\kb/k)})\tors
$$
where $\kb$ is a fixed algebraic closure of $k$, $X_*(T_\kb)_{\Gal(\kb/k)}$
denotes the group of coinvariants, and $(\cdot)\tors$ denotes the 
torsion subgroup.
  If $k$ is a number field and $k_v$ is the completion of $k$ at a place 
$v$, then we define  $C_v(T)=C(T_{k_v})$.
   There is a  canonical map $i_v\colon C_v(T)\to C(T)$ induced by an 
inclusion $\Gal(\kb_v/k_v)\to\Gal(\kb/k)$.
  These definitions were given for connected reductive groups (not only for
tori) by Kottwitz [Ko], see also [BR], 3.4. Kottwitz writes $A(T)$ instead
of $C(T)$.

We compute $C(H)$ for our one-dimensional torus $H$ over $\Q$.
  Clearly 
$$
C(H)=(X_*(H_K)_{\Gal(K/\Q)})\tors=\Z/2\Z\;.
$$
  We have $C_v(H)=1$ if $K\otimes\Q_v$ splits, and $C_v(H)\simeq\Z/2\Z$
if $K\otimes\Q_v$ is a field.
  The map $i_v$ is injective for any $v$.

We now define the local invariants $\kap_v(\ov)$ as in [BR], where
$\ov$ is an orbit of $G(\Q_v)$ in  $X(\Q_v)$.
  The set of orbits of $G(\Q_v)$ in $X(\Q_v)$ is in canonical bijection
with $\ker[H^1(\Q_v,H)\to H^1(\Q_v,G)]$, cf. [Se], I-5.4, 
Cor. 1 of Prop. 36.
  Hence $\ov$ defines a cohomology class $\xi_v\in H^1(\Q_v,H)$.
  The local Tate--Nakayama duality for tori defines a canonical homomorphism
$\bet_v\colon H^1(\Q_v,H)\to C_v(H)$, see [Ko], Thm. 1.2. 
(Kottwitz defines the map $\bet_v$ in a more general setting, when $H$
is any connected reductive group over a number field.)
  The homomorphism $\bet_v$ is an isomorphism for any $v$.
  We set $\kap_v(\ov)=\bet_v(\xi_v)$.
  Note that if $\ov=G(\Q_v)\cdot x^0$, then $\xi_v=0$ and $\kap_v(\ov)=0$;
if $\ov\neq G(\Q_v)\cdot x^0$, then $\xi_v\neq 0$ and $\kap_v(\ov)=1$.

We define the Kottwitz invariant $\kap(\oA)$ of an orbit  $\oA=\prod\ov$
of $G(\A)$ in $X(\A)$ by $\kap(\oA)=\sum_v i_v(\kap_v(\ov))$.
  We identify $C(H)$ with $\Z/2\Z$, and $C_v(H)$ with a subgroup of 
$\Z/2\Z$. 
  With this identifications $\kap(\oA)=\sum\kap_v(\ov)$.

We prefer the multiplicative rather than additive notation.
Instead of $\Z/2\Z$ we consider the group $\{+1,-1\}$, and set
$$
\nu_v(\ov)=(-1)^{\kap_v(\ov)},\ \nu(\oA)=(-1)^{\kap(\oA)}.
$$
 Here $\nu_v(\ov)$ and $\nu(\oA)$ take values $\pm 1$.
 We have $\nu(\oA)=\prod\nu_v(\ov)$.
 Since $\kap_v(\ov)=0$ if and only if $\ov=G(\Q_v)\cdot x^0$, we see 
that $\nu_v(\ov)=+1$ if and only if $\ov=G(\Q_v)\cdot x^0$.   
Hence our $\nu_v(\ov)$ and $\nu(\oA)$ coincide with
$\nu_v(\ov)$ and $\nu(\oA)$ introduced in 1.1.

By Thm. 3.6 of [BR] an adelic orbit $\oA$ contains $\Q$-rational points
if and only if $\kap(\oA)=0$.
  With our multiplicative notation $\kap(\oA)=0$ if and only if 
$\nu(\oA)=+1$.
  Thus $\oA$ contains $\Q$-points if and only if $\nu(\oA)=+1$.
We have deduced Thm. 1.1.1 from [BR], Thm. 3.6.
\qed

\subhead 1.3. Tamagawa measure\endsubhead
We define a gauge form on $X$, i.e. a regular differential form
$\om\in \Lam^2(X)$ without zeroes.
  Recall that $X$ is defined by the equation $f(x)=q$. 
  Choose a differential form $\mu\in\Lam^2 (W)$ such that 
$\mu\wedge df=dx_1\wedge dx_2\wedge dx_3$, where $x_1,x_2,x_3$ are the
coordinates in $W=\Q^3$.
  Let $\om=\mu|_X$, the restriction of $\mu$ to $X$. 
  Then $\om$ is a gauge form on $X$, cf. [BR], 1.3, and it does not depend
on the choice of $\mu$.
  The gauge form $\om$ is $G$-invariant, because there exists a 
$G$-invariant gauge form on $X$, cf. [BR], 1.4,  and a gauge form 
on $X$ is unique up to a scalar multiple, cf. [BR], Cor. 1.5.4.

For any place $v$ of $\Q$ one associates with $\om$ a local measure $m_v$
on $X(\Q_v)$, cf. [We], 2.2.
  We show how to define a Tamagawa measure on $X(\A)$, following 
[BR], 1.6.2.

We have by [BR], 1.8.1, $\mu_p(X)=m_p(X(\Z_p))$, where $\mu_p(X)$ is 
defined in the Introduction. 
  By [We], Thm. 2.2.5, for almost all $p$ we have 
$m_p(X(\Z_p))=\#X(\Fp)$. 

We compute $\#X(\Fp)$. The group $\SO(f)(\Fp)$ acts on $X(\Fp)$ with 
stabilizer $\SO(f')(\F_p)$, where $\SO(f')(\F_p)$ is defined for almost all 
$p$.
  This action is transitive by Witt's theorem.
  Thus $\#X(\Fp)=\#\SO(f)(\Fp)/\#\SO(f')(\Fp)$.
  By [A], III-6, 
$$
\#\SO(f)(\Fp)=p(p^2-1),\quad \#\SO(f')(\Fp)=p-\chi(p),
$$
where $\chi(p)=-1$ if $f'\mod p$ does not represent 0, and 
$\chi(p)=+1$ if $f'\mod p$ represents 0.
  We have $\chi(p)=\dsize\left(\frac{-qD}{p}\right)$.
We obtain for $p\nmid qD$
$$
\#X(\Fp)=\frac{p(p^2-1)}{p-\chi(p)},\quad 
\mu_p(X)=\frac{\#X(\Fp)}{p^2}=\frac{1-{1}/{p^2}}{1-{\chi(p)}/{p}}\;.
$$
  For $p|qD$ set $\chi(p)=0$.
  We define  
$$
L_p(s,\chi)=(1-\chi(p)p^{-s})^{-1},\quad L(s,\chi)=\prod_p L_p(s,\chi)
$$
where $s$ is a complex variable. 
  We set
$$
\lam_p=L_p(1,\chi)^{-1}=1-\frac{\chi(p)}{p},\quad r=L(1,\chi)^{-1}.
$$
  Then the product $\prod_p (\lam_p^{-1}\mu_p)$ converges absolutely, hence 
the family $(\lam_p)$ is a family of convergence factors in the sense of
[We], 2.3.
  We define, as in [BR], 1.6.2, the measures
$$
m_f=r^{-1}\prod_p (\lam_p^{-1}m_p),\quad m=m_\infty m_f\;,
$$
then $m_f$ is a measure on $X(\A_f)$ (where $\A_f$ is the ring of finite
ad\`eles) and $m$ is a measure on $X(\A)$.
  We call $m$ the Tamagawa measure on $X(\A)$.

\subhead 1.4. Counting integer points \endsubhead
  For $T>0$ set $X(\R)^T=\{x\in X(\R): |x|\leq T\}$.

\proclaim{Theorem 1.4.1}
$$
N(T,X)\sim\int_{X(\R)^T\times X(\Zhat)} \del(x)dm.
$$
\endproclaim

In other words,
$$
N(T,X)\sim 2m(\{x\in X(\R)^T\times X(\Zhat): \nu(x)=+1\}).\tag1.4.1
$$

Theorem 1.4.1 follows from [BR], Thm. 5.3 (cf. [BR], 6.4 and [BR], Def. 2.3).

For comparison note that
$$
m(X(\R)^T\times X(\Zhat))=m_\infty(X(\R)^T)m_f(X(\Zhat))
=\mu_\infty(T,X)\SS(X),\tag1.4.2
$$
cf. [BR], 1.8.

The following lemma will be used in the proof of Theorem 0.1.

\proclaim{Lemma 1.4.2} Assume that there exists $y\in X(\R\times\Zhat)$
such that $\nu(y)=+1$. 
  Then the set $X(\Z)$ is infinite.
\endproclaim

\demo{Proof}
Since $\nu$ is a locally constant function on $X(\A)$, there exists
an open subset $\cU_f\in X(\Zhat)$ and an orbit $\cU_\infty$ of $G(\R)$
in $X(\R)$ such that $\nu(x)=+1$ for all $x\in \cU_\infty\times\cU_f$.
  Set $\cU_\infty^T=\{x\in \cU_\infty: |x|\leq T\}$, then 
$m_\infty(\cU_\infty^T)\to\infty$ as $T\to\infty$.
  We have 
$$
\int_{X(\R)^T\times X(\Zhat)} \del(x)dm\geq
\int_{\cU_\infty^T\times\cU_f}\del(x)dm=2m_\infty(\cU_\infty^T)m_f(\cU_f)\;.
$$
Since $2m_\infty(\cU_\infty^T)m_f(\cU_f)\to\infty$ as $T\to \infty$,
we see that 
$$
\int_{X(\R)^T\times X(\Zhat)} \del(x)dm\to\infty \text{ as } T\to\infty,
$$
and by Theorem 1.4.1 $N(T,X)\to\infty$. 
Hence $X(\Z)$ is infinite.
\qed
\enddemo

\subhead 1.5. The constant $c_X$\endsubhead
Here we prove the following result:

\proclaim{Proposition 1.5.1}
$$
N(T,X)\sim c_X \SS(X)\mu_\infty(T,X) \text{ as }T\to\infty
$$
with some constant $c_X$, $0\leq c_X\leq 2$.
\endproclaim

\demo{Proof} If $X(\R)$ has two connected components, then by Theorem 0.2
(which we will prove in Sect. 3 below), $N(T,X)\sim\SS(X)\mu_\infty(T,X)$,
so the proposition holds with $c_X=1$.
 
If $X(\R)$ has one connected component, then $X(\R)$ consists of one
$G(\R)$-orbit and $\nu_\infty(X(\R))=+1$.
  For an orbit $\o_f=\prod\o_p$ of $G(\A_f)$ in $X(\A_f)$ we set 
$\nu_f(\o_f)=\prod_p\nu_p(\o_p)$.
  We regard $\nu_f$ as a locally constant function on $X(\A_f)$ taking
values $\pm 1$.
  We have
$$
\int_{X(\R)^T\times X(\Zhat)} \del(x)dm
=2m_\infty(X(\R)^T)m_f(X(\Zhat)_+)
$$
where $X(\Zhat)_+=\{x_f\in X(\Zhat): \nu_f(x_f)=+1\}$.
Set $c_X=2m_f(X(\Zhat)_+)/m_f(X(\Zhat))$, then $0\leq c_X\leq 2$ and
$$
\int_{X(\R)^T\times X(\Zhat)} \del(x)dm=c_X m_\infty(X(\R)^T) m_f(X(\Zhat))
=c_X\mu_\infty(T,X)\SS(X).
$$
  Using Theorem 1.4.1, we see that
$$
N(T,X)\sim c_X\mu_\infty(T,X)\SS(X) \text{ as } T\to \infty.
$$
\qed
\enddemo

\head 2. An example of $c_X=0$\endhead

Let 
$$
f_1(x_1,x_2, x_3)=-9x_1^2+2x_1x_2+7x_2^2+2x_3^2, \ q=1.
$$
This example was mentioned by Siegel and later mentioned in [BR], 6.4.1.
Here we provide a detailed exposition.

Consider the variety $X$ defined by the equation $f_1(x)=q$.
  We have $f_1(-\half, \half,1)$ $=1$. 
It follows that $f_1$ represents 1 over $\R$ and over $\Z_p$ for $p>2$.
  
We have $f_1(4,1,1)=-127\equiv 1\pmod{2^7}$. 
  We prove that $f_1$ represents 1 over $\Z_2$.
Define a polynomial of one variable  $F(Y)=f_1(4,1,Y)-1,\ F\in \Z_2[Y]$.
Then $F(1)=-2^7$, $|F(1)|_2=2^{-7}$, $F'(Y)=4Y$, $|F'(1)^2|_2=2^{-4}$,
$|F(1)|_2<|F'(1)^2|_2$.
  By Hensel's lemma (cf. [La], II-\S2, Prop. 2) $F$ has a root in $\Z_2$.
  Thus $f_1$ represents 1 over $\Z_2$.

Now we prove that $f_1$ does not represent 1 over $\Z$. 
  I know the following elementary proof from D. Zagier.

We prove the assertion by  contradiction. 
Assume on the contrary that 
$$
-9x_1^2+2x_1x_2+7x_2^2+2x_3^2=1 \text{ for some } x_1, x_2, x_3 \in\Z.
$$
We may write this equation as follows:
$$
2x_3^2-1=(x_1-x_2)^2+8(x_1-x_2)(x_1+x_2).
$$
  The left hand side is odd, hence $x_1-x_2$ is odd and therefore 
$x_1+x_2$ is odd.
  We have $(x_1-x_2)^2\equiv 1 \pmod 8$.
Hence the right hand side is congruent to $1\pmod 8$.
  We see that $x_3$ is odd, hence $2x_3^2-1\equiv 1\pmod{16}$.
But
$$
8(x_1-x_2)(x_1+x_2)\equiv 8 \pmod{16}. 
$$
It follows that 
$$
\gather
(x_1-x_2)^2\equiv 9 \pmod{16}\\
x_1-x_2\equiv\pm 3\pmod8.
\endgather
$$
Therefore $x_1-x_2$ must have a prime factor $p\equiv \pm3 \pmod8$.
Hence $2x_3^2-1$ has a prime factor $p\equiv \pm3 \pmod8$.
On the other hand, if $p|(2x_3^2-1)$, then
$$
2x_3^2\equiv 1\pmod p
$$
and 2 is a square  modulo $p$, $\left(\dsize\frac{2}{p}\right)=1$.
  By the quadratic reciprocity law $p\equiv \pm1\pmod8$. 
  Contradiction.
  We have proved that $f_1$ does not represent 1 over $\Z$,
hence $N(T, X)=0$ for all $T$.

On the other hand, 
$$
\SS(X)\mu_\infty(T, X)=m_f(X(\Zhat))m_\infty(X(\R)^T).
$$
Since $X(\Zhat)$ is a non-empty open subset in $X(\A_f)$, 
$m_f(X(\Zhat))>0$.
  Now the measure $m_\infty(X(R)^T)\to\infty$ as $T\to\infty$.
Hence $\SS(X)\mu_\infty(T,X)\to\infty$ as $T\to\infty$,
and thus $c_{X}=0$.

\head 3. Proof of Theorem 0.1\endhead

\proclaim{Lemma 3.1}
Let $k$ be a field, $\text{\rm char}(k)\neq 2$, and  let $V$ be
a finite-dimensional vector space over $k$.
Let $f$ be a non-degenerate quadratic form on $V$.
Let $u\in \GL(V)(k)$, $f'=u^*f$. Then the map $y\mapsto uy\colon V\to V$
takes the orbits of $\Spin(f)(k)$ in $V$ to the orbits of $\Spin(f')(k)$.
\endproclaim

\demo{Proof}
  Let $x\in V$, $f(x)\neq 0$. 
  The reflection (symmetry) $r_x=r_{f,x}\colon V\to V$
is defined by
$$
r_x(y)=y-\frac{2B(x,y)}{f(x)}x,\quad y\in V,
$$
where $B$ is the symmetric bilinear form on $V$ associated with $f$.
  Every $s\in\SO(f)(k)$ can be written as 
$$
s=r_{x_1}\cdots r_{x_l}\tag3.1
$$
cf. [OM], Thm. 43:3.
\def\SN{{\theta}}
\def\Op{\Theta}
  The spinor norm $\SN(s)$ of $s$ is defined by 
$$
\SN(s)=f(x_1)\cdots f(x_l) \pmod{k^{*2}}\ \in k^*/k^{*2}
$$
and it does not depend on the choice of the representation (3.1),
cf. [OM], \S 55.
Let $\Op(f)$ denote the image of $\Spin(f)(k)$ in $\SO(f)(k)$.
Then $s\in \SO(f)(k)$ is contained in $\Op(f)$ if and only if  $\SN(s)=1$,
cf. [Se], III-3.2 or [Ca], Ch. 10, Thm. 3.3.

Now let $u,f'$ be as above.
  Then $r_{f', ux}=ur_{f,x}u^{-1}$, $f'(ux)=f(x)$, and so 
$\SN_{f'}(usu^{-1})=\SN_f(s)$.
  We conclude that $u\Op(f) u^{-1}=\Op(f')$ and that the map $y\mapsto uy$
takes the orbits of $\Op(f)$ in $V$ to the orbits of $\Op(f')$. 
\qed
\enddemo

Let $f,f'$ be integral-matrix quadratic forms on $\Z^n$ and assume 
that $f'$ is in the genus of $f$.
  Then there exists $u\in\GL_n(\R\times\Zhat)$ such that 
$f'(x)=f(u^{-1}x)$ for $x\in\A^n$.
  Let $q\in\Z$, $q\neq 0$. 
  Let $X$ denote the affine quadric $f(x)=q$, and $X'$ denote the quadric
$f'(x)=q$.
  
\proclaim{ Lemma 3.2} The map $x\mapsto ux\colon \A^n\to\A^n$
takes $X(\R\times\Zhat)$ to $X'(\R\times\Zhat)$ and takes orbits of
$\Spin(f)(\A)$ in $X(\A)$ to orbits of $\Spin(f')(\A)$ in $X'(\A)$.
\endproclaim

\demo{Proof}
Let $A$ denote the matrix of $f$, and $A'$ denote the matrix of $f'$.
  We have
$$
\gather
(u^{-1})^t A u^{-1}=A'\\
A=u^tA' u\;.
\endgather
$$
The variety $X$ is defined by the equation $x^t A x=q$, and $X'$ 
is defined by $x^t A' x=q$.
One can easily check  that the map $x\mapsto ux$ takes $X(\R\times\Zhat)$
to $X'(\R\times\Zhat)$ and $X(\A)$ to $X'(\A)$.

In order to prove that the map $x\mapsto ux\colon X(\A)\to X'(\A)$ takes 
the orbits of $\Spin(f)(\A)$ to the orbits of $\Spin(f')(\A)$, it suffices
to prove that the map $x\mapsto u_v x\colon X(\Q_v)\to X'(\Q_v)$ takes
the orbits of $\Spin(f)(\Q_v)$ to the orbits of $\Spin(f')(\Q_v)$ for
every $v$, where $u_v$ is the $v$-component of $u$. 
This last assertion follows from Lemma 3.1.
\qed
\enddemo

\proclaim{Proposition 3.3} Let $f'$ and $q$ be as in Theorem 0.1, 
in particular $f'$ represents $q$ over $\Z_v$ for any $v$ (we set 
$\Z_\infty=\R$), but not over $\Z$. 
  Let $X'$ be the quadric defined by $f'(x)=q$.
  Then $X'(\R\times\Zhat)$ is contained in one orbit of $\Spin(f')(\A)$.
\endproclaim

\demo{Proof}
Set $G'=\Spin(f')$.
  We prove that $X'(\Z_v)$ is contained in one orbit of $G'(\Q_v)$ for 
every $v$ by contradiction.
  Assume on the contrary that for some $v$  $X'(\Z_v)$  has nontrivial
intersection with two orbits of $G'(\Q_v)$.
  Then $\nu_v$ takes both values $+1$ and $-1$ on $X'(\Z_v)$.
It follows that $\nu$ takes both values $+1$ and $-1$ on $X'(\R\times\Zhat)$.
 Hence by Lemma 1.4.2 $X'$ has infinitely many $\Z$-points.
 This contradicts to the assumption that $f'$ does not represent $q$ over
$\Z$.
\qed
\enddemo

\demo{Proof of Theorem 0.1}
Let  $u\in\GL_3(\R\times\Zhat)$ be such that $f'(x)=f(u^{-1} x)$.
  Let $X,X'$ be as in the beginning of this section,
in particular $X'$ has no $\Z$-points.
  By Prop. 3.3 $X'(\R\times\Zhat)$ is contained in one orbit of 
$\Spin(f')(\A)$.
  It follows from Lemma 3.2 that $X(\R\times\Zhat)$ is contained in one 
orbit of $\Spin(f)(\A)$.
  Since $f$ represents $q$ over $\Z$, this orbit has $\Q$-rational points,
and $\nu$ equals $+1$ on $X(\R\times\Zhat)$.
  Thus $\del$ equals 2 on $X(\R\times\Zhat)$, and by Formulas (1.4.1) and
(1.4.2)
$N(T,X)\sim 2\SS(X)\mu_\infty(T,X)$.
\qed
\enddemo

\head 4. Proof of Theorem 0.2 \endhead

We prove Theorem 0.2.
  We define an involution $\tau_\infty$ of $X(\R)$ by $\tau_\infty(x)=-x$,
$x\in X(\R)\subset \R^3$.
  Since $f(x)=f(-x)$,  $\tau_\infty$ is well defined, i.e takes $X(\R)$
to itself.
  Since $|-x|=|x|$, $\tau_\infty$ takes $X(\R)^T$ to itself.
  We define an involution $\tau$ of $X(\A)$ by defining $\tau$ as 
$\tau_\infty$ on $X(\R)$ and as 1 on $X(\Q_p)$ for all prime $p$.
Then $\tau$ respects the Tamagawa measure $m$ on $X(\A)$.

By assumption $X(\R)$ has two connected components. 
These are two orbits of $\Spin(f)(\R)$.
The involution $\tau_\infty$ of $X(\R)$ interchanges these two orbits.
Thus we have
$$
\gather
\nu_\infty(\tau_\infty(x_\infty))=-\nu_\infty(x_\infty) \text{ for all }
                                    x_\infty\in X(\R)\\
\nu(\tau(x))=-\nu(x) \text{ for all }x\in X(\A)\tag4.1
\endgather
$$

Let $X(\R)_1$ and $X(\R)_2$ be the two connected components of $X(\R)$.
Set
$$
X(R)_1^T=X(R)_1\cap X(\R)^T,\quad X(R)_2^T=X(R)_2\cap X(\R)^T
$$
Then $\tau$ interchanges $X(\R)_1^T\times X(\Zhat)$ and 
$X(\R)_2^T\times X(\Zhat)$. 
From Formula (4.1) we have 
$$
\int_{X(\R)_1^T\times X(\Zhat)} \nu(x)dm
=-\int_{X(\R)_2^T\times X(\Zhat)} \nu(x)dm,
$$
hence
$$
\int_{X(\R)^T\times X(\Zhat)} \nu(x)dm=0.
$$
Since $\del(x)=\nu(x)+1$, we obtain
$$
\int_{X(\R)^T\times X(\Zhat)} \del(x)dm=\int_{X(\R)^T\times X(\Zhat)} dm
=m(X(\R)^T\times X(\Zhat))=\SS(X)\mu_\infty(T,X).
$$
By Theorem 1.4.1
$$
N(T,X)\sim\int_{X(\R)^T\times X(\Zhat)} \del(x)dm.
$$
Thus $N(T,X)\sim\SS(X)\mu_\infty(T,X)$ as $T\to\infty$, i.e. $c_X=1$.
\qed

\comment
\head 5. Numeric calculations for Examples to Theorems 0.1
and 0.2 \endhead
\subhead 5.1. Calculations in Example 0.1.1 \endsubhead
$f(x)=-x_1^2+64x_2^2+2x_3^2$, $D=-128$, $q=1$, hence $-qD\sim 2$ modulo
squares. 
  We have $\chi(p)=\dsize\left(\frac{2}{p}\right)$.

\noindent $\mu_2=2$ (computer calculation).

\noindent $\mu_p=\dsize\frac{1-p^{-2}}{1-\chi(p)/p}$ for $p>2$

\noindent $\SS(X)=2\dsize\frac{L(1,\chi)}{\frac{3}{4}\zeta(2)}$ and 
$\zeta(2)=\pi^2/6$.

\noindent Take $|x|=(x_1^2+64x_2^2+2x_3^2)^{1/2}$, 
          then $\mu_\infty(T,X)\sim\pi T/8$.

\noindent $L(1,\chi)\approx 0.623$ (computer calculation).

\noindent $2\SS(X)\mu_\infty(T,X)\sim 4 L(1,\chi) T/\pi\approx 0.794 T$.

\noindent Numeric calculation gives $N(T,X)/T=0.8$ for $T=10000$.

\subhead 5.2. Calculations in Example 0.2.1 \endsubhead
Here $f(x)=-x_1^2+x_2^2+x_3^2$, $q=-1$, $D=-1$, $-qD=-1$.
We have
$$
\chi(p)=\left(\frac{-1}{p}\right)=(-1)^{(p-1)/2} \text{ for }p>2.
$$
For any odd $n$ set $\chi(n)=(-1)^{(n-1)/2}$, for even $n$ set $\chi(n)=0$.
  We have:

\noindent $\mu_2=\frac{1}{2}$ (numeric calculation for $k=3$);

\noindent $\mu_p=\dsize\frac{1-p^{-2}}{1-\chi(p)/p}$ for $p>2$;

\noindent  $\SS(X)=\dsize\frac{1}{2}\frac{L(1,\chi)}{\frac{3}{4}\zeta(2)}$;

\noindent $\zeta(2)=\pi^2/6$,  $L(1,\chi)=1-\frac{1}{3}+\frac{1}{5}-\cdots
                                          =\pi/4$,

\noindent so $\SS(X)=1/\pi$.

Take $|x|=(x_1^2+x_2^2+x_3^2)^{1/2}$, then

\noindent $\mu_\infty(T,X)\sim\sqrt{2}\pi T \text{ as }T\to\infty$,

\noindent so
$\SS(X)\mu_\infty(T,X)\sim\sqrt{2}T \text{ as }T\to\infty$.

\noindent Numeric calculation gives $N(T,X)/T=1.4072$ for $T=10000$.

\endcomment



\enddocument